\documentclass[10pt]{article}

\usepackage{amsmath}
\usepackage{amsfonts}
\usepackage{amssymb}
\usepackage[english]{babel}
\usepackage{amsthm}

\def\pf{\par\noindent {\bf Proof}~\par\noindent}
\def\qed{~\hfill{$\square$}\pagebreak[1]\par\medskip\par}

\newcommand{\mR}{\mathbb{R}}
\newcommand{\mC}{\mathbb{C}}
\newcommand{\mN}{\mathbb{N}}

\newcommand{\mS}{\mathbb{S}}

\newcommand{\mcM}{\mathcal{M}}

\newcommand{\gf}{\mathfrak{f}}
\newcommand{\gfd}{\mathfrak{f}^{\dagger}}
\newcommand{\gsl}{\mathfrak{sl}}

\newcommand{\ol}{\overline}
\newcommand{\olz}{\ol{z}}

\newcommand{\uX}{\underline{X}}
\newcommand{\ux}{\underline{x}}

\newcommand{\uz}{\underline{z}}
\newcommand{\uzd}{\underline{z}^{\dagger}}
\newcommand{\tuz}{\widetilde{\uz}}
\newcommand{\tuzd}{\tuz^{\dagger}}

\newcommand{\upz}{\partial_{\uz}}

\newcommand{\upzd}{\partial_{\uz}^{\dagger}}

\newcommand{\nz}{\vert\underline{z}\vert}

\newcommand{\p}{\partial}
\newcommand{\dirac}{\underline{\p}}

\newtheorem{theorem}{Theorem}

\newtheorem{proposition}{Proposition}
\newtheorem{remark}{Remark}
\newtheorem{corollary}{Corollary}

\begin{document}
\title{Embedding Factors for Branching \\ in Hermitian Clifford Analysis}
\author{F.\ Brackx$^\ast$, H.\ De Schepper$^\ast$, R.\ L\'{a}vi\v{c}ka$^\ddagger$ \& V.\ Sou\v{c}ek$^\ddagger$}

\date{\small{$^\ast$ Clifford Research Group, Faculty of Engineering, Ghent University\\
Building S22, Galglaan 2, B-9000 Gent, Belgium\\
$^\ddagger$ Mathematical Institute, Faculty of Mathematics and Physics, Charles University\\
Sokolovsk\'a 83, 186 75 Praha, Czech Republic}}

\maketitle

\begin{abstract}
\noindent A step 2 branching decomposition of spaces of homogeneous Hermitian monogenic polynomials in $\mC^n$ is established with explicit embedding factors in terms of the generalized Jacobi polynomials, which allows for an inductive construction of an orthogonal basis for those spaces. The embedding factors and the orthogonal bases are fully worked out in the complex dimension 2 case, with special interest for the Appell property.
\end{abstract}



\section{Introduction}


This paper deals with the decomposition of spaces of homogeneous Hermitian monogenic polynomials in $\mC^n$, in terms of spaces of  ditto polynomials in one complex dimension less, i.e. in $\mC^{n-1}$. The obtained decomposition formula can be seen as a step 2 branching (in real dimensions). Hermitian monogenic functions form one of the actual research topics in Clifford analysis, which, in its most basic form, is a higher dimensional generalization of holomorphic function theory in the complex plane, and, at the same time, a refinement of harmonic analysis, see e.g. \cite{bds, gilbert, dss, guerleb, ghs}. At the heart of Clifford analysis lies the notion of a monogenic function, i.e.\ a Clifford algebra valued null solution of the Dirac operator $\dirac = \sum_{\alpha=1}^m e_{\alpha} \, \p_{X_{\alpha}}$, where $(e_1, \ldots, e_m)$ is an orthonormal basis of $\mR^m$ underlying the construction of the real Clifford algebra $\mR_{0,m}$. We refer to this setting as the Euclidean one, since the fundamental group leaving the Dirac operator $\dirac$ invariant is the orthogonal group $\mbox{O}(m;\mR)$, which is doubly covered by the Pin($m$) group of the Clifford algebra $\mR_{0,m}$. In case the dimension $m$ is even, say $m=2n$, so--called Hermitian Clifford analysis was recently introduced as a refinement of Euclidean Clifford analysis (see the books \cite{rocha,struppa} and the series of papers \cite{sabadini,eel,partI,partII,toulouse,eel2,howe}). The considered functions now take values in the complex Clifford algebra $\mC_{2n}$ or in complex spinor space $\mS_n$. Hermitian Clifford analysis is based on the introduction of an additional datum, a (pseudo) complex structure $J$, inducing an associated Dirac operator $\dirac_J$; it then focusses on the simultaneous null solutions of both operators $\dirac$ and $\dirac_J$, called Hermitian monogenic functions. The corresponding function theory still is in full development, see also \cite{hehe, hermwav, mabo, mahi, cama, eel3}. It is worth mentioning that the traditional holomorphic functions of several complex variables are a special case of Hermitian monogenic functions taking values in a specific homogeneous part of spinor space $\mS_n$.\\[-2mm]

To meet the needs for numerical calculations, recently much effort has been put into the construction of orthogonal bases for spaces of  homogeneous monogenic polynomials, mostly called spherical monogenics, in the framework of both Euclidean and Hermitian Clifford analysis. Indeed, the basis polynomials, sometimes called Fueter polynomials, appearing in the Taylor series expansion of monogenic functions, are not useful for that purpose since they are not orthogonal with respect to the natural $L_2$--inner product. Explicit constructions of orthogonal polynomial bases in the Euclidean Clifford analysis context were carried out in e.g. \cite{cacao,cagubo} in a direct analytic way starting from spherical harmonics, and in e.g. \cite{bock, lav, lav1, lav2,  lasova,  pvl} by the so--called Gel'fand--Tsetlin [GT] approach. The notion of GT--basis stems from group representation theory: every irreducible finite dimensional module over a classical Lie group has its GT--basis (see e.g. \cite{molev}), the construction of which is based on the branching of the corresponding spaces of spherical monogenics. Note that branching for e.g. spherical monogenics in $\mR^m$ is de facto a direct sum decomposition in products of spherical monogenics in $\mR^{m-1}$ multiplied by certain embedding factors. The GT--approach was also used in \cite{rhovla,rhoro,rhofb,rhohds} for designing orthogonal bases for spaces of homogeneous Hermitian monogenic polynomials, with special attention for  the so--called Appell property of the constructed bases in complex dimension $2$. Basis elements are said to possess the Appell property if their derivatives are, up to a multiplicative constant, again basis elements;  this property is, quite naturally, important for numerical applications. \\[-2mm]

In this paper the final step is made by explicitly determining the embedding factors for the step 2 branching of spaces of homogeneous Hermitian monogenic polynomials. Branching formulae in traditional Clifford analysis were already established in \cite{dss} and in \cite{lasova} for the decomposition of these spaces under the action of Spin($m_1) \times$ Spin$(m_2)$, with $m_1+m_2$ the dimension of the considered Euclidean space, leading to the inductive construction of orthogonal bases.
Starting point in our method precisely is the decomposition formula of spaces of standard spherical monogenics in $\mR^m$ in terms of spherical monogenics in $\mR^{m_1}$ and $\mR^{m_2}$, with $m_1+m_2=m$, which is recalled in Section 3 and adapted to the Hermitian Clifford analysis framework. Next, the Fischer decomposition of spaces of standard spherical monogenics in terms of Hermitian spherical monogenics is used.  Adequate combination of both decomposition formulae lead to the explicit embedding factors. The newly obtained decomposition formula allows for the inductive construction of orthogonal bases for spaces of Hermitian spherical monogenics. The embedding factors and orthogonal basis polynomials are described in detail in the complex dimension 2 case, with special emphasis on the Appell property. To make the paper self--contained an introduction on Clifford analysis is included.


\section{Preliminaries on Clifford analysis}


For a detailed description of the structure of a Clifford algebra we refer to e.g.\ \cite{port}. Here we only recall the necessary basic notions. The real Clifford algebra $\mathbb{R}_{0,m}$ is constructed over the vector space $\mathbb{R}^{0,m}$ endowed with a non--degenerate quadratic form of signature $(0,m)$ and generated by the orthonormal basis $(e_1,\ldots,e_m)$. The non--commutative Clifford or geometric multiplication in $\mathbb{R}_{0,m}$ is governed by the rules 
\begin{equation}\label{multirules}
e_{\alpha} e_{\beta} + e_{\beta} e_{\alpha} = -2 \delta_{\alpha \beta} \ \ , \ \ \alpha,\beta = 1,\ldots ,m
\end{equation}
As a basis for $\mathbb{R}_{0,m}$ one takes for any set $A=\lbrace j_1,\ldots,j_h \rbrace \subset \lbrace 1,\ldots,m \rbrace$, the element $e_A = e_{j_1} \ldots e_{j_h}$, with $1\leq j_1<j_2<\cdots < j_h \leq m$, together with $e_{\emptyset}=1$, the identity element. Any Clifford number $a$ in $\mathbb{R}_{0,m}$ may thus be written as $a = \sum_{A} e_A a_A$, $a_A \in \mathbb{R}$, or still as $a = \sum_{k=0}^m \lbrack a \rbrack_k$, where $\lbrack a \rbrack_k = \sum_{|A|=k} e_A a_A$ is the so--called $k$--vector part of $a$. Euclidean space $\mathbb{R}^{0,m}$ is embedded in $\mathbb{R}_{0,m}$ by identifying $(X_1,\ldots,X_m)$ with the Clifford vector $\uX = \sum_{\alpha=1}^m e_{\alpha}\, X_{\alpha}$, for which it holds that $\uX^2 = - |\uX|^2$. The vector valued first order differential operator $\dirac = \sum_{\alpha=1}^m e_{\alpha}\, \p_{X_{\alpha}}$, called Dirac operator, is the Fourier or Fischer dual of the Clifford variable $\uX$. It is this operator which underlies the notion of monogenicity of a function, a notion which is the higher dimensional counterpart of holomorphy in the complex plane. More explicitly, a function $f$ defined and continuously differentiable in an open region $\Omega$ of $\mathbb{R}^{m}$ and taking values in (a subspace of) the Clifford algebra $\mathbb{R}_{0,m}$, is called (left) monogenic in $\Omega$ if $\dirac\lbrack f \rbrack = 0$ in $\Omega$. As the Dirac operator factorizes the Laplacian: $\Delta_m = - \dirac^2$, monogenicity can be regarded as a refinement of harmonicity. The Dirac operator being rotationally invariant, this framework is usually referred to as Euclidean Clifford analysis.\\[-2mm]

When allowing for complex constants, the generators $(e_1,\ldots, e_{m})$, still satisfying \eqref{multirules}, produce the complex Clifford algebra $\mathbb{C}_{m} = \mathbb{R}_{0,m} \oplus i\, \mathbb{R}_{0,m}$. Any complex Clifford number $\lambda \in \mathbb{C}_{m}$ may thus be written as $\lambda = a + i b$, $a,b \in \mathbb{R}_{0,m}$, leading to the definition of the Hermitian conjugation $\lambda^{\dagger} = (a +i b)^{\dagger} = \overline{a} - i \overline{b}$, where the bar notation stands for the Clifford conjugation in $\mathbb{R}_{0,m}$, i.e. the main anti--involution for which $\overline{e}_{\alpha} = -e_{\alpha}$, $\alpha=1, \ldots,m$. This Hermitian conjugation leads to a Hermitian inner product on $\mathbb{C}_{m}$ given by $(\lambda,\mu) = \lbrack \lambda^{\dagger} \mu \rbrack_0$ and its associated norm $|\lambda| = \sqrt{ \lbrack \lambda^{\dagger} \lambda \rbrack_0} = ( \sum_A |\lambda_A|^2 )^{1/2}$. This is the framework for Hermitian Clifford analysis, which emerges from Euclidean Clifford analysis by introducing a so--called complex structure, i.e.\ an $\mbox{SO}(m;\mR)$--element $J$ with $J^2=-\mathbf{1}$ (see \cite{partI,partII}), forcing the dimension to be even; from now on we put $m=2n$. Usually $J$ is chosen to act upon the generators of $\mC_{2n}$ as $J[e_j] = -e_{n+j}$ and $J[e_{n+j}] = e_j$, $j=1,\ldots,n$. By means of the projection operators $\pm \frac{1}{2}(\mathbf{1} \pm iJ)$ associated to $J$, first the Witt basis elements $(\gf_j,\gf_j^{\dagger})^n_{j=1}$ for $\mC_{2n}$ are obtained: 
\begin{eqnarray*}
\gf_j & = & \phantom{-}\frac{1}{2} (\mathbf{1} + iJ) [e_j] = \phantom{-}\frac{1}{2} (e_j - i \, e_{n+j}), \qquad j=1,\ldots,n \\
\gf_j^\dagger & = & -\frac{1}{2} (\mathbf{1} - iJ)[e_j] =  -\frac{1}{2} (e_j + i\, e_{n+j}), \qquad j=1,\ldots,n
\end{eqnarray*}
The Witt basis elements satisfy the respective Grassmann and duality identities
$$
\gf_j \gf_k + \gf_k \gf_j = \gf_j^{\dagger} \gf_k^{\dagger} + \gf_k^{\dagger} \gf_j^{\dagger} = 0, \quad
\gf_j \gf_k^{\dagger} + \gf_k^{\dagger} \gf_j = \delta_{jk}, \ \ j,k=1,\ldots, n
$$
whence they are isotropic: $(\gf_j)^2 = 0, (\gfd_j)^2 = 0, j=0,\ldots,n$.
Next, denoting a vector in $\mR^{0,2n}$ by $(x_1, \ldots , x_n, y_1 , \ldots , y_n)$, which is  identified with the Clifford vector $\uX  =  \sum_{j=1}^n (e_j\, x_j + e_{n+j}\, y_j)$, the Hermitian Clifford variables $\uz$ and $\uzd$ are produced similarly:
$$
\uz = \frac{1}{2} (\mathbf{1} + iJ) [\uX] = \sum_{j=1}^n \gf_j\, z_j, \quad
\uzd = -\frac{1}{2} (\mathbf{1} - iJ)[\uX] = \sum_{j=1}^n \gf_j^{\dagger}\, \olz_j
$$
where complex variables $z_j = x_j + i y_j$ have been introduced, with complex conjugates $\olz_j = x_j - i y_j$, $j=1,\ldots,n$. Finally, the Euclidean Dirac operator $\dirac$ gives rise to the Hermitian Dirac operators $\upz$ and $\upzd$:
$$
\upzd = \frac{1}{4} (\mathbf{1} + iJ) [\dirac] = \sum_{j=1}^n \gf_j\, \p_{\olz_j}, \quad
\upz = -\frac{1}{4} (\mathbf{1} - iJ)[\dirac] = \sum_{j=1}^n \gf_j^{\dagger}\, \p_{z_j}  
$$
involving the Cauchy--Riemann operators $\p_{\olz_j} = \frac{1}{2} (\p_{x_j} + i \p_{y_j})$  and their complex conjugates $\p_{z_j} = \frac{1}{2} (\p_{x_j} - i \p_{y_j})$  in the $z_j$--planes, $j=1,\ldots,n$. Observe that Hermitian vector variables and Dirac operators are isotropic, i.e. $\uz^2 = (\uzd)^2 = 0$ and $(\upz)^2 = (\upzd)^2 = 0$, whence the Laplacian allows for the decomposition and factorization
$$
\Delta_{2n} = 4(\upz \upzd + \upzd \upz) = 4(\upz + \upzd)^2 = -4(\upzd - \upz)^2
$$
while dually
\begin{displaymath}
-(\uz-\uzd)^2 = (\uz + \uzd)^2 = \uz\, \uzd + \uzd \uz = |\uz|^2 = |\uzd|^2 = |\uX|^2
\end{displaymath}

\vspace*{3mm}
We consider functions with values in an irreducible representation $\mS_n$ of $\mC_{2n}$, called spinor space, which is realized within $\mC_{2n}$ using a primitive idempotent $I = I_1 \ldots I_n$, with $I_j = \gf_j \gf_j^{\dagger}$, $j=1,\ldots,n$. With that choice, $\gf_jI=0, j=1,\ldots,n$, and so $\mS_n \equiv \mC_{2n} I \cong {\bigwedge}_n^{\dagger} I$, where ${\bigwedge}_n^{\dagger} = \bigwedge(\gfd_1, \ldots, \gfd_n)$ denotes the Grassmann algebra generated by $(\gfd_1, \ldots, \gfd_n)$. Hence $\mS_n$ decomposes as 
$$
\mS_n = \bigoplus\limits_{r=0}^n \ \mS_n^{(r)} = \bigoplus\limits_{r=0}^n \ ( {\bigwedge}_n^{\dagger} )^{(r)} I
$$
with $( {\bigwedge}_n^\dagger )^{(r)} = \mbox{span}_{\mC} ( \gfd_{k_1} \wedge \gfd_{k_2} \wedge \cdots \wedge \gfd_{k_r} : \{ k_1,\ldots,k_r \} \subset \{1,\ldots,n  \} )$.
By singling out one of the Witt basis vectors, viz. $\gfd_n$, and still using the same idempotent $I$, we can consider $\mS_{n-1} = \bigwedge(\gfd_1, \ldots, \gfd_{n-1})I$, which leads to the direct sum decomposition of spinor space
$$
\mS_n = \mS_{n-1} \oplus \gfd_n \, \mS_{n-1}
$$
and of its homogeneous parts
$$
\mS_n^{(r)} = ({\bigwedge}_{n-1}^{\dagger})^{(r)}(\gfd_1, \ldots, \gfd_{n-1}) I \oplus \gfd_n \ ({\bigwedge}_{n-1}^{\dagger})^{(r-1)}(\gfd_1, \ldots, \gfd_{n-1}) I
$$

A continuously differentiable function $g$ in an open region $\Omega$ of $\mathbb{R}^{2n}$ with values in (a subspace of) $\mathbb{C}_{2n}$ then is called (left) Hermitian monogenic (or h--monogenic) in $\Omega$ if and only if it satisfies in $\Omega$ the system $\upz\, g = 0 = \upzd\, g$, or, equivalently, the system $\dirac \, g = \dirac_J \, g$, with $\dirac_J = J[\dirac]$. A major difference between Hermitian and Euclidean Clifford analysis concerns the underlying group invariance, which for $(\upz,\upzd)$ breaks down to the group U$(n)$, see e.g.\ \cite{partI,partII}. This plays a fundamental role in the construction of orthogonal bases for spaces of Hermitian monogenic polynomials.\\[-2mm]

The space of homogeneous monogenic polynomials on $\mC^{n}$, taking values in spinor space $\mS_n$ and with global degree of homogeneity $k$ in the complex variables $(\uz,\uzd)$, is denoted by $\mcM_k(n)$. When specifying the bidegree $(a,b)$ in $(\uz,\uzd)$ we denote the corresponding space by $\mcM_{a,b}(n)$, while $\mcM_{a,b}^{r}(n)$ stands for its subspace where the values are taken in the homogeneous part $\mS_n^{(r)}$. 

\section{The embedding factors}

Let us consider the above defined space $\mcM_k(n)$. It is well known that this space may be decomposed into U$(n)$--irreducibles involving appropriate spaces $\mcM_{a,b}^{r}(n)$. This is one of the so--called Fischer decompositions, for which we refer to \cite{rhoro} --where it is proved via analytic methods-- and to \cite{damdeef} for a group representation approach. This Fischer decomposition explicitly reads as follows.

\begin{theorem}
The space $\mcM_k(n)$ may be decomposed into {\rm U}($n$)-irreducibles as
\begin{equation}
\label{monoghermonog-decomp}
\mcM_k(n) =  \bigoplus_{a+b=k} \, \bigoplus_{r=0}^n \ \mcM_{a,b}^{r}(n) \ \oplus \ \bigoplus_{a+b=k-1} \, \bigoplus_{r=1}^{n-1} \, \left( (b+n-r)\uz+ (a+r)\uzd \right) \mcM_{a,b}^{r}(n)
\end{equation}
\end{theorem}

There is however a second decomposition for this space via step 2 branching (in real dimension), involving spaces of homogeneous monogenic polynomials defined in $\mR^{2n-2} \cong \mC^{n-1}$. Starting point is the formula established in \cite{dss}, where orthogonal bases for the spaces $\mcM_k(\mR^{m}; \mR_{0,m})$ of monogenic $k$--homogeneous polynomials defined in Euclidean space $\mR^m$ and taking values in the real Clifford algebra $\mR_{0,m}$, are constructed in an inductive way based on the splitting $\mR^{m} = \mR^{m_1} \times \mR^{m_2}$, with $m_1 + m_2 =m$. We recall the expressions for the embedding factors established in \cite[Lemma 4.4]{dss}, pp. 260-262, with some minor corrections. For $i,j,\ell\in\mN$ we have that
\begin{eqnarray}
S_{2\ell,j,i}(\ux_1,\ux_2)&=&(|\ux_1|^2+|\ux_2|^2)^{\ell-1}\left((|\ux_1|^2+|\ux_2|^2)P_\ell^{i+\frac{m_2}{2}-1,j+\frac{m_1}{2}-1}\left(\frac{|\ux_1|^2-|\ux_2|^2}{|\ux_1|^2+|\ux_2|^2}\right)\right . \nonumber \\
&& \left . -\ux_1\ux_2P_{\ell-1}^{i+\frac{m_2}{2},j+\frac{m_1}{2}}\left(\frac{|\ux_1|^2-|\ux_2|^2}{|\ux_1|^2+|\ux_2|^2}\right)\right) \label{Seven}
\end{eqnarray}
and
\begin{eqnarray}
S_{2\ell+1,j,i}(\ux_1,\ux_2)&=&(|\ux_1|^2+|\ux_2|^2)^{\ell}\left((\ell+j+\frac{m_1}{2})\ux_2P_\ell^{i+\frac{m_2}{2},j+\frac{m_1}{2}-1}\left(\frac{|\ux_1|^2-|\ux_2|^2}{|\ux_1|^2+|\ux_2^2}\right)\right . \nonumber \\
&& \left .  -(\ell+i+\frac{m_2}{2})\ux_1P_{\ell}^{i+\frac{m_2}{2}-1,j+\frac{m_1}{2}}\left(\frac{|\ux_1|^2-|\ux_2|^2}{|\ux_1|^2+|\ux_2|^2}\right)\right) \label{Soneven}
\end{eqnarray} 
where $\ux_1 \in \mR^{m_1}$, $\ux_2 \in \mR^{m_2}$, with $m_1+m+2 = m$, and $P_\ell^{\alpha,\beta}$ denote the classical Jacobi polynomials in one real variable given by
$$
P_\ell^{\alpha,\beta}(t) = \sum_{s=0}^{\ell} \binom{\ell + \alpha}{s} \binom{\ell + \beta}{\ell - s} \left( \frac{t-1}{2}  \right)^{\ell-s}  \left( \frac{t+1}{2}   \right)^{s}
$$

\begin{remark}
The above expressions (\ref{Seven}) and (\ref{Soneven}) do not literally correspond to the ones in \cite{dss}, since the notation for the indices has been changed in order to facilitate further use. In particular, the new first index now indicates the total degree of the polynomial.  
\end{remark}

Now taking $m=2n, m_1=2, m_2=2n-2$, introducing the complex Clifford variables $\uz = \tuz + \gf_n z_n$, $\uzd = \tuzd + \gfd_n  \olz_n$ and restricting the function values to the respective spinor spaces $\mS_n$ and $\mS_{n-1}$ for which holds the direct sum decomposition $\mS_n = \mS_{n-1} \oplus \gfd_n \, \mS_{n-1}$ explained in the foregoing section, the following decomposition theorem is obtained.

\begin{theorem}
The space $\mcM_k(n)$ of \, $\mS_n$--valued spherical monogenics of degree $k$ may be decomposed as
\begin{equation}
\label{monogmonog-decomp}
\mcM_k(n) =   \bigoplus_{0 \leq i+j \leq k}  S_{k-j-i,j,i}(\uz, \uzd)  (\olz_n)^j  \, \mcM_{i}(n-1) \oplus  S_{k-j-i,j,i}(\uz, \uzd) (z_n)^j  \gfd_n  \, \mcM_{i}(n-1) 
\end{equation}
where the polynomials $S_{k-j-i,j,i}$ appearing in the embedding factors are given by
\begin{equation}
\label{Seven}
S_{2\ell,j,i}(\uz, \uzd) = |\uz|^{2\ell-2} \left(  |\uz|^2 P_{\ell}^{i+n-2,j}(t) - (\gf_n z_n - \gfd_n \olz_n)(\tuz-\tuzd)P_{\ell-1}^{i+n-1,j+1}(t) \right)
\end{equation}
and
\begin{equation}
\label{Sodd}
S_{2\ell+1,j,i}(\uz, \uzd) = |\uz|^{2\ell} \left((\ell+j+1)(\tuz - \tuzd)P_{\ell}^{i+n-1,j}(t)-(\ell+i+n-1) (\gf_n z_n - \gfd_n \olz_n)P_{\ell}^{i+n-2,j+1}(t)  \right)
\end{equation}
with $t$ representing the dimensionless variable
$$
t = \frac{|z_n|^2 - |\tuz|^2}{|\uz|^2}
$$
\end{theorem}

To give a more concrete idea what this decomposition \eqref{monogmonog-decomp} looks like, we establish its explicit form for the cases $k=0, 1, 2$.
For $k=0$ we obtain, quite trivially,
$$
\mcM_0(n) =  \mcM_{0}(n-1) \oplus \gfd_n  \, \mcM_{0}(n-1) 
$$
For $k=1$ we find
\begin{eqnarray*}
\mcM_1(n) &=& \left(\tuz - \tuzd - (n-1)(\gf_n z_n -\gfd_n \olz_n)\right) \mcM_{0}(n-1) \\
&& \oplus \, \left(\tuz - \tuzd - (n-1)(\gf_n z_n -\gfd_n \olz_n)\right) \gfd_n \mcM_{0}(n-1)\\
&& \oplus \, \mcM_{1}(n-1) \, \oplus \, \gfd_n \mcM_{1}(n-1) \, \oplus \, \olz_n  \mcM_{0}(n-1) \, \oplus \, z_n  \gfd_n  \mcM_{0}(n-1)
\end{eqnarray*}
For $k=2$ this decomposition already counts $12$ terms:
\begin{eqnarray*}
\mcM_{2}(n) &=& \left((n-1)z_n\olz_n - |\tuz|^2 - (\gf_n z_n -\gfd_n \olz_n)(\tuz - \tuzd)\right) \, \mcM_0(n-1) \\
&& \oplus \left((n-1)z_n\olz_n - |\tuz|^2 - (\gf_n z_n -\gfd_n \olz_n)(\tuz - \tuzd)\right) \, \gfd_n \, \mcM_0(n-1) \\
&& \oplus \left(\tuz - \tuzd - n \, (\gf_n z_n -\gfd_n \olz_n)\right) \,\mcM_1(n-1) \oplus \left(\tuz - \tuzd - n \, (\gf_n z_n -\gfd_n \olz_n)\right) \, \gfd_n \, \mcM_1(n-1) \\
&& \oplus \left(2(\tuz - \tuzd) - (n-1) (\gf_n z_n -\gfd_n \olz_n)\right) \,\olz_n \, \mcM_0(n-1) \\
&& \oplus \left(2(\tuz - \tuzd) - (n-1) (\gf_n z_n -\gfd_n \olz_n)\right) \, z_n \, \gfd_n \, \mcM_0(n-1) \\
&& \oplus \ \mcM_2(n-1) \oplus \gfd_n \, \mcM_2(n-1) \oplus \ \olz_n \, \mcM_1(n-1) \oplus z_n \, \gfd_n \, \mcM_1(n-1) \\
&& \oplus \ \olz_n^2 \, \mcM_0(n-1) \oplus z_n^2 \, \gfd_n \, \mcM_0(n-1)
\end{eqnarray*}
It may be checked by direct, yet far from trivial, computation, that each of the components in the above direct sum decomposition indeed is monogenic w.r.t.\ the Dirac operator $\dirac = 2(\upzd - \upz)$.\\

Now let us have a look at the decomposition aimed at. We want to decompose the space  $\mcM_{a,b}^{r}(n)$ of $(a,b)$--homogeneous Hermitian monogenic polynomials in the variables $(\uz, \uzd)$ and with values in $\mS_n^{( r )}$, in terms of spaces of homogeneous Hermitean monogenic polynomials in the variables $(\tuz, \tuzd)$. From the start we assume that $0<r<n$, the cases $r=0$ and $r=n$ being treated separately (see Remark 3). It follows from representation theory and from the Cauchy-Kovalevskaya extension in Hermitian Clifford analysis, see \cite{CK,rhovla}, that this decomposition must have the following form:
\begin{equation}
\label{step2-decomp}
\mcM_{a,b}^{r}(n) =  \bigoplus_{c=0}^a \, \bigoplus_{d=0}^b \  X_{a,b;c,d}^{r,r} \, \mcM_{c,d}^{r} (n-1)  \oplus
                                 \bigoplus_{c=0}^a \, \bigoplus_{d=0}^b \  X_{a,b;c,d}^{r,r-1} \, \mcM_{c,d}^{r-1} (n-1)
\end{equation}
but the problem remains to explicitly determine the embedding factors $X_{a,b;c,d}^{r,s}$, $c=0,\ldots,a$, $d=0,\ldots,b$, $s=r,r-1$. When applying the Fischer decomposition \eqref{monoghermonog-decomp}  on each of the spaces $\mcM_{i}(n-1)$ appearing in the decomposition \eqref{monogmonog-decomp}, we find
\begin{eqnarray}
\mcM_k(n) & = & \bigoplus_{0 \leq i+j \leq k}  S_{k-j-i,j,i}(\uz, \uzd)  (\olz_n)^j  \, (
  \bigoplus_{c+d=i} \, \bigoplus_{r=0}^{n-1} \ \mcM_{c,d}^{r}(n-1)  \nonumber \\
& & \hspace*{2cm} \oplus \  \bigoplus_{c+d=i-1} \, \bigoplus_{r=1}^{n-2} \, \left( (d+n-1-r)\uz + (c+r)\uzd \right) \mcM_{c,d}^{r}(n-1)  ) \nonumber\\
& &   \bigoplus_{0 \leq i+j \leq k}  S_{k-j-i,j,i}(\uz, \uzd) (z_n)^j  \gfd_n  \, (
   \bigoplus_{c+d=i}\, \bigoplus_{r=0}^{n-1} \ \mcM_{c,d}^{r}(n-1)  \nonumber \\
 & &   \hspace*{2cm}  \oplus \  \bigoplus_{c+d=i-1} \, \bigoplus_{r=1}^{n-2} \, \left( (d+n-1-r)\uz +  (c+r)\uzd \right) \mcM_{c,d}^{r}(n-1)  ) 
 \label{monogmonog-decomp-bis}
\end{eqnarray}

When comparing the decompositions \eqref{monoghermonog-decomp} and \eqref{monogmonog-decomp-bis}, it becomes clear that each embedding factor $X_{a,b;c,d}^{r,s}$ can be expressed as an appropriate combination of two embedding factors appearing in \eqref{monogmonog-decomp-bis}. Indeed, we can only combine irreducible pieces of the same representation character with respect to U$(n-1) \ \times$ U$(1)$, that is with the same labels $r,c,d$ and with the same factor $z_n^j$ or $\olz_n^j$. Putting $a-c=u$ and $b-d=v$, we have to distinguish between the cases $u < v$, $u=v$ and $u > v$.\\[-2mm]

\noindent {\bf The case $u < v$}\\
We put $u = \ell$ and $v = \ell +j$, or $a = c + \ell$ and $b = d + \ell +j$, and we aim at determining the embedding factor $X_{a,b;c,d}^{r,r}$ which maps 
$\mcM_{c,d}^{r}(n-1)$ into $\mcM_{c+\ell,d+\ell+j}^{r}(n)$. Thence its action should increase the degree in the $\uz$ variables by $\ell$, increase the degree in the $\uzd$ variables by $\ell+j$, meanwhile keeping the homogeneity degree of spinor space unaltered. To achieve this we propose to act on $\mcM_{c,d}^{r}(n-1)$ with an expression of the form
$$
X_{c+\ell,d+\ell+j;c,d}^{r,r} = \alpha_1 \, S_{2\ell,j,c+d} (\olz_n)^j + \beta_1 \, S_{2\ell-1,j,c+d+1} (\olz_n)^j ((d+n-1-r)\tuz + (c+r)\tuzd)
$$
where $\alpha_1$ and $\beta_1$ are real scalars to be determined in such a way that the {\em bad} terms appearing in this expression, i.e. the terms which do not respect the bidegree and homogeneity degree aimed at, vanish. Using the explicit forms \eqref{Seven}--\eqref{Sodd} of the $S$--polynomials, the resulting equation reads:
$$
- \alpha_1 |\uz|^{2\ell-2} (\olz_n)^{j+1} P_{l-1}^{c+d+n-1,j+1} \gfd_n  \tuzd + \beta_1 |\uz|^{2\ell-2} (\olz_n)^{j+1} (\ell+c+d+n-1) P_{l-1}^{c+d+n-1,j+1} (c+r) \gfd_n \tuzd
$$
from which it follows that, up to constants:
$$
\alpha_1 = (\ell+c+d+n-1)(c+r), \quad \beta_1 = 1
$$
Note that the above reasoning and the obtained result remain valid in the case where $j=0$ or $u=v$. For the embedding factor $X_{a,b;c,d}^{r,r-1}$, mapping $\mcM_{c,d}^{r-1}(n-1)$ to $\mcM_{c+\ell,d+\ell+j}^{r}(n)$, we propose, following a similar reasoning:
$$
X_{c+\ell,d+\ell+j;c,d}^{r,r-1} = \alpha_2 \, S_{2\ell,j-1,c+d+1} (\olz_n)^{j-1} ((d+n-r)\tuz + (c+r-1)\tuzd) + \beta_2 \, S_{2\ell+1,j-1,c+d} (\olz_n)^{j-1} 
$$
Again substituting the explicit forms \eqref{Seven}\eqref{Sodd} for the $S$--polynomials, this expression for the desired embedding factor contains two bad terms, the vanishing of which leads to
$$
\alpha_2 = \ell+j, \quad \beta_2 = -(d+n-r)
$$
Note that the obtained result is only valid when $j>0$, so the case $u=v$ is excluded here.\\

\noindent {\bf The case $u > v$}\\
Now we put $v = \ell$ and $u = \ell + j$, or $a = c + \ell + j$ and $b = d + \ell$. In a similar way as above we find that the embedding factor $X_{a,b;c,d}^{r,r}$ mapping $\mcM_{c,d}^{r}(n-1)$ to $\mcM_{c+\ell+j,d+\ell}^{r}(n)$ is given by
$$
X_{c+\ell+j,d+\ell;c,d}^{r,r} = \alpha_3 \, S_{2\ell,j-1,c+d+1} (z_n)^{j-1} \gfd_n ((d+n-1-r)\tuz + (c+r)\tuzd) + \beta_3 \, S_{2\ell+1,j-1,c+d} (z_n)^{j-1} \gfd_n
$$
with
$$
\alpha_3 = \ell+j, \quad \beta_3 = -(c+r)
$$
and this result is only valid when $j>0$, meaning that the strict inequality $u>v$ should be respected. The embedding factor $X_{a,b;c,d}^{r,r-1}$, mapping $\mcM_{c,d}^{r-1}(n-1)$ to $\mcM_{c+\ell+j,d+\ell}^{r}(n)$, is found to be
$$
X_{c+\ell+j,d+\ell;c,d}^{r,r-1} = \alpha_4 \, S_{2\ell,j,c+d} (z_n)^j \gfd_n + \beta_4 \, S_{2\ell-1,j,c+d+1} (z_n)^j \gfd_n ((d+n-r)\tuz + (c+r-1)\tuzd)
$$
with
$$
\alpha_4 = (\ell+c+d+n-1)(d+n-r), \quad \beta_4 = 1
$$
a result which remains valid for $j=0$ or $u=v$.\\

In this way we have proved the following branching theorem.

\begin{theorem}
The U$(n)$--module $\mcM_{a,b}^{r}(n)$  of polynomials in the variables $(\uz,\uzd)$ may be decomposed into  U$(n-1)$--irreducibles of homogeneous Hermitian monogenic polynomials in the variables $(\tuz,\tuzd)$ as
\begin{equation}
\label{step2-decompfinal}
\mcM_{a,b}^{r}(n) = 
\bigoplus_{c=0}^a \, \bigoplus_{d=0}^b \  X_{a,b;c,d}^{r,r} \, \mcM_{c,d}^{r} (n-1)  \oplus \bigoplus_{c=0}^a \, \bigoplus_{d=0}^b \  X_{a,b;c,d}^{r,r-1} \, \mcM_{c,d}^{r-1} (n-1) 
\end{equation}
the embedding factors being given by\\[1mm]
(i) for $a-c \leq b-d, a=c+\ell, b=d+\ell+j$
\begin{eqnarray*}
&& \hspace*{-5mm} X_{c+\ell,d+\ell+j;c,d}^{r,r} \\
&& = (\ell+n+c+d-1)(c+r)  S_{2\ell,j,c+d} (\olz_n)^j +  S_{2\ell-1,j,c+d+1} (\olz_n)^j ((d+n-1-r)\tuz + (c+r)\tuzd)
\end{eqnarray*}
(ii) for $a-c < b-d, a=c+\ell, b=d+\ell+j \ (j \geq 1)$
\begin{eqnarray*}
&& \hspace*{-5mm}  X_{c+\ell,d+\ell+j;c,d}^{r,r-1} \\
&& =  (\ell+j)  S_{2\ell,j-1,c+d+1} (\olz_n)^{j-1} ((d+n-r)\tuz + (c+r-1)\tuzd)  -   (d+n-r)  S_{2\ell+1,j-1,c+d} (\olz_n)^{j-1} 
\end{eqnarray*}
(iii) for $a-c > b-d, a = c+\ell+j \ (j \geq 1), b= d+\ell$
\begin{eqnarray*}
&& \hspace*{-5mm}  X_{c+\ell+j,d+\ell;c,d}^{r,r} \\
&& = (\ell+j)  S_{2\ell,j-1,c+d+1} (z_n)^{j-1} \gfd_n ((d+n-1-r)\tuz + (c+r)\tuzd)  -  (c+r)  S_{2\ell+1,j-1,c+d} (z_n)^{j-1} \gfd_n
\end{eqnarray*}
(iv) for $a-c \geq b-d, a = c+\ell+j, b= d+\ell$
\begin{eqnarray*}
&& \hspace*{-5mm} X_{c+\ell+j,d+\ell;c,d}^{r,r-1} \\
&&   =  (\ell+n+c+d-1)(d+n-r)  S_{2\ell,j,c+d} (z_n)^j \gfd_n +  S_{2\ell-1,j,c+d+1} (z_n)^j \gfd_n ((d+n-r)\tuz + (c+r-1)\tuzd)
\end{eqnarray*}
\end{theorem}

For the lowest dimensional cases, we will now give the explicit form of this decomposition. For $a=0, b=0$ and $0 < r < n$, we find
\begin{eqnarray*}
X_{0,0;0,0}^{r,r} &=& (n-1) r \, S_{0,0,0} \\
X_{0,0;0,0}^{r,r-1} &=& (n-1)(n-r) S_{0,0,0} \gfd_n
\end{eqnarray*}
which leads to the trivially expected decomposition
$$
\mcM_{0,0}^r(n) = \mcM_{0,0}^r(n-1) \oplus \gfd_n \mcM_{0,0}^{r-1}(n-1)
$$
For $a=1, b=0$ and $0 < r < n$, we find
\begin{eqnarray*}
X_{1,0;0,0}^{r,r} &=&  S_{0,0,1} \, \gfd_n \, ((n-1-r)\tuz + r \tuzd) - r \, S_{1,0,0} \, \gfd_n = P_0^{n-1,0} (n-1) \, \gfd_n \, \tuz + r(n-1) P_0^{n-2,1} z_n \gf_n \gfd_n \\
X_{1,0;0,0}^{r,r-1} &=& (n-1)(n-r) S_{0,1,0} \, z_n \, \gfd_n = (n-1)(n-r)P_0^{n-2,1} z_n \, \gfd_n\\
X_{1,0;1,0}^{r,r} &=& n (r+1) S_{0,0,1} = n (r+1)P_0^{n-1,0} \\
X_{1,0;1,0}^{r,r-1} &=& n (n-r) S_{0,0,1} \, \gfd_n = n (n-r)P_0^{n-1,0}\, \gfd_n
\end{eqnarray*}
leading to the decomposition
$$
\mcM_{1,0}^r(n) =  \mcM_{1,0}^r(n-1) \, \oplus \, \gfd_n \mcM_{1,0}^{r-1}(n-1) \, \oplus  \, (\gfd_n \tuz + r z_n) \mcM_{0,0}^r(n-1) \, \oplus  \, z_n \gfd_n \mcM_{0,0}^{r-1}(n-1)
$$
It can be verified that the four components indeed are spaces of Hermitian monogenic polynomials.\\

\noindent In a similar way we obtain for $a=0, b=1$ and $0 < r < n$, 
\begin{eqnarray*}
X_{0,1;0,0}^{r,r} &=& (n-1) r \, S_{0,1,0} \olz_n = (n-1) r \, P_0^{n-2,1} \olz_n \\
X_{0,1;0,0}^{r,r-1} &=& S_{0,0,1} ((n-r)\tuz + (r-1)\tuzd) -(n-r)S_{1,0,0}  \\
& = & (n-1) P_0^{n-1,0} \tuzd + (n-r)(n-1) P_0^{n-2,1}(\gf_n z_n - \gfd_n \olz_n) \\
X_{0,1;0,1}^{r,r} &=& n \, r \, S_{0,0,1} = n \, r \, P_0^{n-1,0} \\
X_{0,1;0,1}^{r,r-1} &=& n(n-r+1) S_{0,0,1} \gfd_n = n(n-r+1) P_0^{n-1,0} \gfd_n
\end{eqnarray*}
leading to the decomposition
$$
\mcM_{0,1}^r(n) =  \mcM_{0,1}^r(n-1) \, \oplus \, \gfd_n \mcM_{0,1}^{r-1}(n-1) \,
\oplus  \, \olz_n \, \mcM_{0,0}^r(n-1) \, \oplus  \, ( \tuzd - (n-r) \olz_n \gfd_n) \mcM_{0,0}^{r-1}(n-1)
$$

\begin{remark}
The decomposition formula \eqref{step2-decompfinal} can be used for constructing, by induction on the dimension, an orthogonal basis of the space $\mcM_{a,b}^r(n)$, called Gel'fand--Tsetlin basis. This induction process takes off in the complex plane ($n=1$) where the orthogonal bases are well known: both spaces $\mcM_{0,b}^0(1)$ and $\mcM_{a,0}^1(1)$
are one--dimensional and spanned by $\frac{(\olz_1)^b}{b!}I$ and $\frac{(z_1)^a}{a!}\gfd_1I$ respectively. The basis polynomials constructed in this way, will be expressed in terms of the classical Jacobi polynomials. This approach should be compared with the Cauchy--Kovalevskaya procedure, used in \cite{rhohds}, where the basis polynomials are expressed as natural powers of $z_1, \ldots, z_n$ and $\olz_1, \ldots, \olz_n$.
\end{remark}

\begin{remark}
In the above considerations it was assumed, from the beginning, that $0<r<n$. When $r=0$, Hermitian monogenicity is nothing but anti--holomorphy and the polynomials at stake only depend on the variables $(\olz_1, \ldots, \olz_n)$, while for $r=n$ the Hermitian monogenic polynomials are holomorphic and only depend on the variables $(z_1, \ldots, z_n)$ (see e.g. \cite{partII}). In these two exceptional cases the branching decomposition formula \eqref{step2-decompfinal} takes a specific form:
\begin{eqnarray*}
\mcM_{0,b}^0(n) &=& \bigoplus_{d=0}^b \, X_{0,b;0,d}^{0,0} \, \mcM_{0,d}^0(n-1) \quad {\rm with} \quad X_{0,b;0,d}^{0,0} = (\olz_n)^{b-d} \\
\mcM_{a,0}^n(n) &=& \bigoplus_{c=0}^a \, X_{a,0;c,0}^{n,n-1} \, \mcM_{c,0}^{n-1}(n-1) \quad {\rm with} \quad X_{a,0;c,0}^{n,n-1} = (z_n)^{a-c} \, \gfd_n
\end{eqnarray*}
\end{remark}


\section{The case of complex dimension 2}

In this section we will decompose the spaces $\mcM_{a,b}^{(r)}(\mC^2)$ of spherical Hermitian monogenics in the complex variables $z_1, z_2, \olz_1, \olz_2$, taking values in the homogeneous parts $\mS_2^{(r)}$ of spinor space $\mS_2$, according to the branching formula established in the previous section. In this way it is possible to construct at once an orthogonal basis for these spaces  $\mcM_{a,b}^{(r)}(\mC^2)$, since the bases in one complex dimension lower, i.e.\ in the complex plane, are trivially well--known. Needless to say that the orthogonal bases are important for applications in real dimension $4$. Note that they were already explicitly constructed in \cite{rhohds} by using the so--called Gel'fand--Tsetlin approach, but here they will be expressed in terms of the classical Jacobi polynomials.

\subsection{The embedding factors in complex dimension 2}

As already mentioned, the induction procedure for the construction of an orthogonal basis starts in the complex plane, with the well--known bases given in Remark 2 for the one--dimensional spaces $\mcM_{0,b}^0(1)$ and $\mcM_{a,0}^1(1)$.\\

In complex dimension 2, spinor space $\mS_2$ decomposes into three homogeneous parts 
$$
\mS_2^{(0)} = {\rm span}_{\mC}\{1\}I, \quad \mS_2^{(1)} = {\rm span}_{\mC}\{\gfd_1, \gfd_2 \}I \quad \mbox{and} \quad \mS_2^{(2)} = {\rm span}_{\mC}\{\gfd_1 \gfd_2  \}I
$$
with $I = \gf_1\gfd_1\gf_2\gfd_2$. In the case where $r=0$ we obtain the branching
$$
\mcM_{0,b}^0(2) = \bigoplus_{d=0}^b \ (\olz_2)^{b-d} \ \mcM_{0,d}^0(1)
$$
It can be readily checked that both sides show dimension $b+1$. Moreover the right hand side generates the following orthogonal basis for $\mcM_{0,b}^0(2)$:
$\left\{(\olz_2)^{b-d} (\olz_1)^d I: d = 0, 1, \ldots, b\right\}$. In the case where $r=2$ we obtain the branching
$$
\mcM_{a,0}^2(2) = \bigoplus_{c=0}^a \ (z_2)^{a-c}\, \gfd_2 \ \mcM_{c,0}^1(1)
$$
both sides clearly showing dimension $a+1$. Moreover the right hand side generates the following orthogonal basis for $\mcM_{a,0}^2(2)$: $
\left\{(z_2)^{a-c} (z_1)^c \, \gfd_1 \gfd_2 : c = 0, 1, \ldots, a\right\}$.  In the case where $r=1$ we obtain the branching
\begin{equation}
\label{n=2,r=1}
\mcM_{a,b}^1(2) = \bigoplus_{c=0}^a \   X_{a,b;c,0}^{1,1}   \ \mcM_{c,0}^1(1)  \oplus \bigoplus_{d=0}^b \   X_{a,b;0,d}^{1,0}   \ \mcM_{0,d}^0(1)
\end{equation}
where the dimension of both sides is clearly seen to be $a+b+2$. Let us first have a look at the embedding factor $X_{a,b;c,0}^{1,1}$. In the case where $a \leq b$, there holds $u=a-c \leq b-0 =v$
and so
$$
X_{a,b;c,0}^{1,1} = (a+2)(c+1) \, S_{2(a-c),b-a+c,c} \, \olz_2^{b-a+c} + S_{2(a-c)-1,b-a+c,c+1} \, \olz_2^{b-a+c} \, (c+1) \, \tuzd
$$
When acting on the space $\mcM_{c,0}^1(1)$ the second term has no contribution and hence, up to constants,
\begin{eqnarray*}
&& \hspace*{-7mm} X_{a,b;c,0}^{1,1} \, \mcM_{c,0}^1(1) \\
&& = \olz_2^{b-a+c} \left( \nz^{2(a-c)}  P_{a-c}^{c,b-a+c}(t)  - \nz^{2(a-c-1)}  (\gf_2 z_2 - \gfd_2 \olz_2)(\gf_1 z_1 - \gfd_1 \olz_1)  P_{a-c-1}^{c+1,b-a+c+1}(t)  \right)  \mcM_{c,0}^1(1)  \\
&& = \left( \olz_2^{b-a+c}  \nz^{2(a-c)} P_{a-c}^{c,b-a+c}(t) +  z_1 \olz_2^{b-a+c+1}  \nz^{2(a-c-1)} P_{a-c-1}^{c+1,b-a+c+1}(t)  \gfd_2  \gf_1 \right)  \mcM_{c,0}^1(1) 
\end{eqnarray*}
with
$$
t = \frac{|z_2|^2 - |z_1|^2}{|\uz|^2} =  \frac{z_2 \olz_2 - z_1 \olz_1}{|\uz|^2} 
$$
In the case where $a >b$ the first term in the branching \eqref{n=2,r=1} splits into
$$
\bigoplus_{c=0}^{a-b-1} \   X_{a,b;c,0}^{1,1}   \, \mcM_{c,0}^1(1)  \oplus \bigoplus_{c=a-b}^a \   X_{a,b;c,0}^{1,1}  \, \mcM_{c,0}^1(1)
$$
For the second part we can use the above expression for the embedding factor $X_{a,b;c,0}^{1,1} $, while for the first part we obtain
$$
X_{a,b;c,0}^{1,1} = (a-c) S_{2b,a-c-b-1,c+1} z_2^{a-c-b-1} \gfd_2 (c+1) \tuzd  -  (c+1) S_{2b+1,a-c-b-1,c} z_2^{a-c-b-1} \gfd_2
$$
When this embedding factor acts on the space $\mcM_{c,0}^1(1)$ then the first term has no contribution, leading, up to constants, to
\begin{eqnarray*}
&& \hspace*{-7mm} X_{a,b;c,0}^{1,1} \,\mcM_{c,0}^1(1) \\
&& \hspace*{-7mm} =  \nz^{2b}  z_2^{a-c-b-1} \left(  (a-c)(\tuz-\tuzd) P_b^{c+1,a-c-b-1}(t)
- (b+c+1)(\gf_2 z_2 - \gfd_2 \olz_2 P_b^{c,a-c-b}(t)     \right)  \gfd_2  \mcM_{c,0}^1(1) \\
&& \hspace*{-7mm} = \left( (a-c) z_1 z_2^{a-c-b-1} \nz^{2b} P_b^{c+1,a-c-b-1}(t) \gf_1 \gfd_2 - (b+c+1) z_2^{a-c-b} \nz^{2b} P_b^{c,a-c-b}(t) \right)  \mcM_{c,0}^1(1)
\end{eqnarray*}
Let us now turn our attention to the embedding factor $X_{a,b;0,d}^{1,0}$. If $a \geq b$ then $u = a \geq b-d = v$ and 
$$
X_{a,b;0,d}^{1,0} = (b+1)(d+1) S_{2(b-d),a-b+d,d} z_2^{a-b+d} \gfd_2  + S_{2(b-d)-1,a-b+d,d+1} z_2^{a-b+d} \gfd_2 (d+1) \tuz
$$
yielding, up to constants,
\begin{eqnarray*}
&& \hspace*{-10mm}  X_{a,b;0,d}^{1,0}  \, \mcM_{0,d}^0(1) \\
&& \hspace*{-7mm} = \left( \nz^{2(b-d)} P_{b-d}^{d,a-b+d}(t)  -  \nz^{2(b-d-1)} (\gf_2 z_2 - \gfd_2 \olz_2)(\tuz - \tuzd) P_{b-d-1}^{d+1,a-b+d+1}(t)  \right)  \, z_2^{a-b+d} \gfd_2 \mcM_{0,d}^0(1) \\
&&  \hspace*{-7mm} = \left( \nz^{2(b-d)} z_2^{a-b+d} P_{b-d}^{d,a-b+d}(t) \gfd_2 - \nz^{2(b-d-1)}  \olz_1 z_2^{a-b+d+1} P_{b-d-1}^{d+1,a-b+d+1}(t) \gfd_1
\right)  \, \mcM_{0,d}^0(1)
\end{eqnarray*}
In the case where $a < b$ we have to distinguish between $d$ running from $0$ till $b-a-1$ and hence $a < b-d$, and $d$ running from $b-a$ till $b$ and hence $a \geq b-d$.
For $d=b-a, \ldots, b$ we can use the above expression for $X_{a,b;0,d}^{1,0} $, whereas for $d=0, \ldots, b-a-1$ we obtain
$$
X_{a,b;0,d}^{1,0} = (b-d) S_{2a, b-d-a-1,d+1} \olz_2^{b-d-a-1} (d+1)\tuz - (d+1) S_{2a+1,b-d-a-1,d} \olz_2^{b-d-a-1}
$$
and hence, up to constants,
\begin{eqnarray*}
&& \hspace*{-7mm}  X_{a,b;0,d}^{1,0}  \, \mcM_{0,d}^0(1) \\
&& \hspace*{-7mm} = \nz^{2a} \left( (b-d) (\tuz - \tuzd) P_a^{d+1,b-d-a-1}(t) - (a+d+1)(\gf_2 z_2 - \gfd_2 \olz_2) P_a^{d,b-d-a}(t) \right) \olz_2^{b-d-a-1}   \mcM_{0,d}^0(1) \\
&& \hspace*{-7mm}  =  \left( (b-d) \nz^{2a} \olz_1 \olz_2^{b-d-a-1} P_a^{d+1,b-d-a-1}(t) \gfd_1 - (a+d+1) \nz^{2a} \olz_2^{b-d-a} P_a^{d,b-d-a}(t) \gfd_2 \right)   \mcM_{0,d}^0(1)
\end{eqnarray*}
In conclusion, the branching \eqref{n=2,r=1} for the space $\mcM_{a,b}^1(2)$ takes the following form:\\

\noindent
(i) if $a<b$ then
$$
\mcM_{a,b}^1(2) = \bigoplus_{c=0}^a \   X_{a,b;c,0}^{1,1}   \ \mcM_{c,0}^1(1)  \oplus \bigoplus_{d=0}^{b-a-1} \   X_{a,b;0,d}^{1,0}   \ \mcM_{0,d}^0(1) \oplus \bigoplus_{d=b-a}^{b} \   X_{a,b;0,d}^{1,0}   \ \mcM_{0,d}^0(1) 
$$
or, with explicit embedding factors,
\begin{eqnarray*}
&& \hspace*{-7mm}  \mcM_{a,b}^1(2) = \\
&& \hspace*{-2mm}  \bigoplus_{c=0}^a \left( \olz_2^{b-a+c}  \nz^{2(a-c)} P_{a-c}^{c,b-a+c}(t) +  z_1 \olz_2^{b-a+c+1}  \nz^{2(a-c-1)} P_{a-c-1}^{c+1,b-a+c+1}(t)  \gfd_2  \gf_1 \right)  \mcM_{c,0}^1(1) \\
&& \hspace*{-7mm} \oplus \bigoplus_{d=0}^{b-a-1}  \left( (b-d) \nz^{2a} \olz_1 \olz_2^{b-d-a-1} P_a^{d+1,b-d-a-1}(t) \gfd_1 
- (a+d+1) \nz^{2a} \olz_2^{b-d-a} P_a^{d,b-d-a}(t) \gfd_2 \right)   \mcM_{0,d}^0(1) \\
&& \hspace*{-7mm} \oplus \bigoplus_{d=b-a}^{b} \left( \nz^{2(b-d)} z_2^{a-b+d} P_{b-d}^{d,a-b+d}(t) \gfd_2 - \nz^{2(b-d-1)}  \olz_1 z_2^{a-b+d+1} P_{b-d-1}^{d+1,a-b+d+1}(t) \gfd_1
\right)  \mcM_{0,d}^0(1)
\end{eqnarray*}

\noindent
(ii) if $a>b$ then
$$
\mcM_{a,b}^1(2) = \bigoplus_{c=0}^{a-b-1} \   X_{a,b;c,0}^{1,1}   \, \mcM_{c,0}^1(1)  \oplus \bigoplus_{c=a-b}^a \   X_{a,b;c,0}^{1,1}  \, \mcM_{c,0}^1(1) \oplus \bigoplus_{d=0}^b \   X_{a,b;0,d}^{1,0}   \ \mcM_{0,d}^0(1)
$$
or, with explicit embedding factors,
\begin{eqnarray*}
&& \hspace*{-7mm} \mcM_{a,b}^1(2) = \\
&& \hspace*{-4mm} \bigoplus_{c=0}^{a-b-1}  \left( (a-c) z_1 z_2^{a-c-b-1} \nz^{2b} P_b^{c+1,a-c-b-1}(t) \gf_1 \gfd_2 - (b+c+1) z_2^{a-c-b} \nz^{2b} P_b^{c,a-c-b}(t) \right)  \mcM_{c,0}^1(1) \\
&& \hspace*{-7mm} \oplus \bigoplus_{c=a-b}^a  \left( \olz_2^{b-a+c}  \nz^{2(a-c)} P_{a-c}^{c,b-a+c}(t) + z_1 \olz_2^{b-a+c+1}  \nz^{2(a-c-1)} P_{a-c-1}^{c+1,b-a+c+1}(t)  \gfd_2  \gf_1 \right)  \mcM_{c,0}^1(1)  \\
&& \hspace*{-7mm} \oplus \, \, \bigoplus_{d=0}^b \ \left( \nz^{2(b-d)} z_2^{a-b+d} P_{b-d}^{d,a-b+d}(t) \gfd_2 - \nz^{2(b-d-1)}  \olz_1 z_2^{a-b+d+1} P_{b-d-1}^{d+1,a-b+d+1}(t) \gfd_1 \right)   \mcM_{0,d}^0(1)
\end{eqnarray*}

\noindent
(iii) if $a=b$ then
\begin{eqnarray*}
\mcM_{a,a}^1(2) &=& \bigoplus_{c=0}^{a} \   X_{a,a;c,0}^{1,1}   \, \mcM_{c,0}^1(1) \oplus \bigoplus_{d=0}^a \   X_{a,a;0,d}^{1,0}   \, \mcM_{0,d}^0(1) \\
&=& \bigoplus_{c=0}^{a} \ \left( \nz^{2(a-c)} \olz_2^c P_{a-c}^{c,c}(t) +  \nz^{2(a-c-1)} \olz_2^{c+1} z_1 P_{a-c-1}^{c+1,c+1}(t) \gfd_2 \gf_1 \right) \, \mcM_{c,0}^1(1) \\
& &  \oplus \ \bigoplus_{d=0}^a \ \left(\nz^{2(a-d)} z_2^d P_{a-d}^{d,d}(t) \gfd_2  -  \nz^{2(a-d-1)} z_2^{d+1} \olz_1 P_{a-d-1}^{d+1,d+1}(t) \gfd_1
\right) \, \mcM_{0,d}^0(1)
\end{eqnarray*}

\subsection{Orthogonal basis in complex dimension 2}

Using the known orthogonal bases for the spaces $\mcM_{c,0}^1(1)$ and $\mcM_{0,d}^0(1)$, the above branching formulae yield the following orthogonal basis for $\mcM_{a,b}^1(2) $:\\[-2mm]

\noindent (i) if $a<b$ then
\begin{eqnarray*}
&& \hspace*{-6mm} \mcM_{a,b}^1(2) = \\[-1mm]
&& \hspace*{-3mm}{\rm span}_{c=0 \ldots a} \  \left( z_1^c \olz_2^{b-a+c}  \nz^{2(a-c)} P_{a-c}^{c,b-a+c}(t) \gfd_1 I + z_1^{c+1} \olz_2^{b-a+c+1}  \nz^{2(a-c-1)} P_{a-c-1}^{c+1,b-a+c+1}(t)  \gfd_2 I \right) \\
&& \hspace*{-3mm} \oplus  {\rm span}_{d=0 \ldots b-a-1}  \left( (b-d) \nz^{2a}  \olz_1^{d+1} \olz_2^{b-d-a-1} P_a^{d+1,b-d-a-1}(t) \gfd_1 I \right .\\[-1mm]
&& \hspace*{50mm} \left . - (a+d+1) \nz^{2a} \olz_1^d \olz_2^{b-d-a} P_a^{d,b-d-a}(t) \gfd_2 I \right) \\
&& \hspace*{-3mm} \oplus  {\rm span}_{d=b-a \ldots b}  \left( \nz^{2(b-d)}  \olz_1^d z_2^{a-b+d} P_{b-d}^{d,a-b+d}(t) \gfd_2 I - \nz^{2(b-d-1)}  \olz_1^{d+1}  z_2^{a-b+d+1} P_{b-d-1}^{d+1,a-b+d+1}(t) \gfd_1 I \right) 
\end{eqnarray*}

\noindent
(ii) if $a>b$ then
\begin{eqnarray*}
&& \hspace*{-6mm} \mcM_{a,b}^1(2) = \\[-1mm]
&& \hspace*{-3mm}{\rm span}_{c=0 \ldots a-b-1}   \left( (a-c) z_1^{c+1}  z_2^{a-c-b-1} \nz^{2b} P_b^{c+1,a-c-b-1}(t)  \gfd_2 I \right .\\[-1mm]
&& \hspace*{50mm} \left . + (b+c+1)  z_1^c  z_2^{a-c-b} \nz^{2b} P_b^{c,a-c-b}(t) \gfd_1 I \right) \\
&&  \hspace*{-3mm}\oplus  {\rm span}_{c=a-b \ldots a}  \left(  z_1^c  \olz_2^{b-a+c}  \nz^{2(a-c)} P_{a-c}^{c,b-a+c}(t) \gfd_1 I  +  z_1^{c+1}  \olz_2^{b-a+c+1}  \nz^{2(a-c-1)} P_{a-c-1}^{c+1,b-a+c+1}(t)  \gfd_2  I  \right) \\
&& \hspace*{-3mm}\oplus  {\rm span}_{d=0 \ldots b} \left( \nz^{2(b-d)}  \olz_1^d  z_2^{a-b+d} P_{b-d}^{d,a-b+d}(t) \gfd_2 I - \nz^{2(b-d-1)}  \olz_1^{d+1}  z_2^{a-b+d+1} P_{b-d-1}^{d+1,a-b+d+1}(t) \gfd_1 I \right) 
\end{eqnarray*}

\noindent
(iii) if $a=b$ then
\begin{eqnarray*}
\mcM_{a,a}^1(2)  &=& \bigoplus_{c=0}^{a} \ \left( \nz^{2(a-c)}  z_1^c \olz_2^c P_{a-c}^{c,c}(t) \gfd_1 I +  \nz^{2(a-c-1)}  z_1^{c+1} \olz_2^{c+1} P_{a-c-1}^{c+1,c+1}(t) \gfd_2  I \right)   \\
&& \oplus \ \bigoplus_{d=0}^a \ \left(\nz^{2(a-d)} \olz_1^d z_2^d P_{a-d}^{d,d}(t) \gfd_2  I -  \nz^{2(a-d-1)} z_2^{d+1} \olz_1^{d+1} P_{a-d-1}^{d+1,d+1}(t) \gfd_1  I \right)  
\end{eqnarray*}

We can make these basis polynomials still more explicit by substituting the defining expression for the dimensionless parameter $t$. As for the Jacobi polynomials it holds that 
$$
P_l^{\alpha,\beta}\left(\frac{u}{v}\right) = \frac{1}{v^l} \ \sum_{s=0}^l \ \binom{l+\alpha}{s} \binom{l+\beta}{l-s} \left(\frac{u-v}{2}\right)^{l-s} \left({\frac{u+v}{2}}\right)^s
$$
we obtain
$$
P_l^{\alpha,\beta}(t) = \frac{1}{|\uz|^{2l}} \ \sum_{s=0}^l \ \binom{l+\alpha}{s} \binom{l+\beta}{l-s} \left(-z_1 \ol{z}_1 \right)^{l-s} \left(z_2\ol{z}_2 \right)^s
$$
leading to the introduction of modified Jacobi polynomials in the complex variables $z_1, \ol{z}_1, z_2, \ol{z}_2$
$$
Q_l^{\alpha, \beta}(\uz, \uzd) = |\uz|^{2l} \, P_l^{\alpha,\beta}(t) = \sum_{s=0}^l \ \binom{l+\alpha}{s} \binom{l+\beta}{l-s} \left(-z_1 \ol{z}_1\right)^{l-s} \left(z_2\ol{z}_2\right)^s
$$
These polynomials enjoy the following properties which will be used in the next subsection, and which are proven by direct calculation.

\begin{proposition}
For the polynomials $Q_l^{\alpha, \beta}(\uz, \uzd)$ one has
\begin{itemize}
\item[(i)] $\partial_{z_1} Q_{l}^{\alpha,\beta} = - (l+\beta) \, \ol{z}_1 \, Q_{l-1}^{\alpha+1,\beta}$
\item[(ii)] $\partial_{\ol{z}_1} Q_{l}^{\alpha,\beta} = - (l+\beta) \, z_1 \, Q_{l-1}^{\alpha+1,\beta}$
\item[(iii)] $\partial_{z_2} Q_{l}^{\alpha,\beta} =  (l+\alpha) \, \ol{z}_2 \, Q_{l-1}^{\alpha,\beta+1}$
\item[(iv)] $\partial_{\ol{z}_2} Q_{l}^{\alpha,\beta} =  (l+\alpha) \, z_2 \, Q_{l-1}^{\alpha,\beta+1}$
\item[(v)] $\beta \, Q_l^{\alpha,\beta} + (l+\alpha) \, z_2 \olz_2 \, Q_{l-1}^{\alpha,\beta+1} = (l+\beta) \, Q_l^{\alpha,\beta-1}$
\item[(vi)] $\alpha \, Q_l^{\alpha,\beta} + (l+\beta) \, (- z_1 \olz_1) \, Q_{l-1}^{\alpha+1,\beta} = (l+\alpha) \, Q_l^{\alpha-1,\beta}$
\end{itemize}
\end{proposition}

As an aside notice that appropriate combination of the above properties (v) and (vi) leads to the following recurrence relations for the standard Jacobi polynomials which, as such, we did not encounter in the literature.

\begin{corollary}
The Jacobi polynomials satisfy the following recurrence relations
\begin{itemize}
\item[(i)] $P_{l-1}^{\alpha,\beta} = - \frac{l+\alpha-1}{l+\beta} \, P_{l}^{\alpha-2,\beta}  +  \frac{\alpha-1}{l+\beta} \, P_{l}^{\alpha-1,\beta}  - \frac{\beta-1}{l+\alpha} \, P_{l}^{\alpha,\beta-1} +  \frac{l+\beta-1}{l+\alpha} \, P_{l}^{\alpha,\beta-2}$
\item[(ii)] $ t \, P_{l-1}^{\alpha,\beta} =  \phantom{.} \frac{l+\alpha-1}{l+\beta} \, P_{l}^{\alpha-2,\beta}  -  \frac{\alpha-1}{l+\beta} \, P_{l}^{\alpha-1,\beta}  - \frac{\beta-1}{l+\alpha} \, P_{l}^{\alpha,\beta-1}  +  \frac{l+\beta-1}{l+\alpha} \, P_{l}^{\alpha,\beta-2}$
\item[(iii)] $\frac{1+t}{2} \, P_{l-1}^{\alpha,\beta} =    - \frac{\beta-1}{l+\alpha} \, P_{l}^{\alpha,\beta-1}  +  \frac{l+\beta-1}{l+\alpha} \, P_{l}^{\alpha,\beta-2}$
\item[(iv)] $\frac{1-t}{2} \, P_{l-1}^{\alpha,\beta} = - \frac{l+\alpha-1}{l+\beta} \, P_{l}^{\alpha-2,\beta}  +  \frac{\alpha-1}{l+\beta} \, P_{l}^{\alpha-1,\beta}$
\end{itemize}
\end{corollary}

In terms of the newly introduced polynomials $Q_l^{\alpha, \beta}(\uz, \uzd)$ we may now rewrite the orthogonal basis polynomials for $\mcM_{a,b}^1(2)$ as follows.
\begin{corollary}
The space $\mcM_{a,b}^1(2)$ of spherical Hermitian monogenics in complex dimension 2 shows the following basis:\\[-2mm]

\noindent (i) if $a<b$
\begin{eqnarray*}
p_{a,b;c,0} &=& \frac{1}{a!} \frac{1}{b!} \, \left( z_1^c  \olz_2^{b-a+c}  Q_{a-c}^{c,b-a+c}  \gfd_1  I + 
z_1^{c+1}  \olz_2^{b-a+c+1}  Q_{a-c-1}^{c+1,b-a+c+1}  \gfd_2  I \right), \quad \; \;  c = 0, \ldots, a \\
\widetilde{q}_{a,b;0,d} &=&  \frac{1}{(a+d+1)!} \frac{1}{(b-d)!} \,
\left( (b-d)  \olz_1^{d+1}  \olz_2^{b-d-a-1}  Q_{a}^{d+1,b-d-a-1}  \gfd_1  I \right .\\
&& \hspace*{41mm} \left . -(a+d+1) \olz_1^{d}  \olz_2^{b-d-a}  Q_{a}^{d,b-d-a}  \gfd_2  I \right), \quad d = 0, \ldots, b-a-1\\
q_{a,b;0,d} &=&  \frac{1}{a!} \frac{1}{b!} \, \left( \olz_1^{d+1}  z_2^{a-b+d+1}  Q_{b-d-1}^{d+1,a-b+d+1}  \gfd_1  I - \olz_1^{d}  z_2^{a-b+d}  Q_{b-d}^{d,a-b+d}  \gfd_2  I \right), \quad d = b-a, \ldots, b
\end{eqnarray*}

\noindent (ii) if $a>b$
\begin{eqnarray*}
\widetilde{p}_{a,b;c,0}  &=&  \frac{1}{(a-c)!} \frac{1}{(b+c+1)!} \,
\left( (b+c+1) z_1^{c} z_2^{a-c-b}  Q_{b}^{c,a-c-b}  \gfd_1  I \right . \\
&& \hspace*{33mm} \left . + (a-c) z_1^{c+1} z_2^{a-c-b-1} Q_{b}^{c+1,a-c-b-1} \gfd_2  I \right), \quad c = 0, \ldots, a-b-1\\
p_{a,b;c,0} &=&  \frac{1}{a!} \frac{1}{b!} \, \left( z_1^{c}  \olz_2^{b-a+c}  Q_{a-c}^{c,b-a+c}  \gfd_1  I +
 z_1^{c+1}  \olz_2^{b-a+c+1}  Q_{a-c-1}^{c+1,b-a+c+1}  \gfd_2  I \right), \quad \, c = a-b, \ldots, a\\
q_{a,b;0,d} &=& \frac{1}{a!} \frac{1}{b!} \, \left(  \olz_1^{d+1}  z_2^{a-b+d+1}  Q_{b-d-1}^{d+1,a-b+d+1}  \gfd_1  I - \olz_1^{d}  z_2^{a-b+d}  Q_{b-d}^{d,a-b+d}  \gfd_2  I \right), \; \; \, d = 0, \ldots, b
\end{eqnarray*}

\noindent (iii) if $a=b$
\begin{eqnarray*}
p_{a,a;c,0} &=& \frac{1}{a!} \frac{1}{a!} \, \left( z_1^c \, \olz_2^{c} \, Q_{a-c}^{c,c} \, \gfd_1 \, I + 
z_1^{c+1} \, \olz_2^{c+1} \, Q_{a-c-1}^{c+1,c+1} \, \gfd_2 \, I \right), \quad \; \, c = 0, \ldots, a\\
q_{a,a;0,d} &=&  \frac{1}{a!} \frac{1}{a!} \, \left( \olz_1^{d+1} \, z_2^{d+1} \, Q_{a-d-1}^{d+1,d+1} \, \gfd_1 \, I - \olz_1^{d} \, z_2^{d} \, Q_{a-d}^{d,d} \, \gfd_2 \, I \right), \quad d = 0, \ldots, a
\end{eqnarray*}
\end{corollary}

\subsection{The Appell property in complex dimension 2}

In this subsection we show that in complex dimension $n=2$, the above constructed orthogonal bases of spherical Hermitian monogenics possess the Appell property with respect to all variables, that is, by differentiating any basis polynomial with respect to one of the variables $\olz_2$, $z_2$, $\olz_1$ or $z_1$, always a multiple of another basis element is obtained; it is even so that by a suitable choice of normalizing factors no multiplicative constants are needed. This property is obvious for the $\mS_2^{(0)}$-- and $\mS_2^{(2)}$--valued basis polynomials, while for the $\mS_2^{(1)}$--valued polynomials an explicit calculation will be carried out. In fact the Appell property for the derivatives with respect to $\olz_2$ and $z_2$ holds in any dimension $n$, since differentiation with respect to the ''last variables''  $\overline{z}_n$ and $z_n$ is obviously $U(n-1)$-invariant. So the following result confirms the Appell property for the derivatives with respect to $\olz_2$ and $z_2$, and proves the Appell property for the derivatives with respect to $\olz_1$ and $z_1$. The significance of the Appell property is the following. Considering a finite dimensional subspace of spherical Hermitian monogenics with bidegree of homogeneity bounded by fixed constants $a$ and $b$, each of the four derivatives is represented with respect to the orthogonal basis by a very simple nilpotent matrix: it is a block matrix (with respect to the decomposition given by the irreducible pieces of the considered subspace) where almost all blocks are zero matrices and where each non--zero block has the property that every column contains at most one nontrivial entry. In this way the Appell property clearly makes numerical calculations very efficient.

\begin{proposition}
For the derivatives with respect to $z_2$ and $\olz_2$ of the orthogonal basis polynomials of the space $\mcM_{a,b}^1(2)$ of spherical Hermitian monogenics in complex dimension 2, there holds:\\[-2mm]

\noindent (i) if $a<b$\\[-2mm]

\noindent
* $\partial_{z_2} p_{a,b;c,0} = p_{a-1,b;c,0}, \quad c = 0, \ldots, a-1$\\
* $\partial_{z_2} p_{a,b;a,0} = 0$\\
* $\partial_{z_2} \widetilde{q}_{a,b;0,d} = \widetilde{q}_{a-1,b;0,d}, \quad d = 0, \ldots, b-a-1$\\
* $\partial_{z_2} q_{a,b;0,b-a} = \widetilde{q}_{a-1,b;0,b-a}$ \\
* $\partial_{z_2} q_{a,b;0,d} = q_{a-1,b;0,d}, \quad d = b-a+1, \ldots, b$\\[-2mm]

\noindent
* $\partial_{\olz_2} p_{a,b;c,0} = p_{a,b-1;c,0}, \quad c = 0, \ldots, a$\\
* $\partial_{\olz_2} \widetilde{q}_{a,b;0,d} = \widetilde{q}_{a,b-1;0,d}, \quad d = 0, \ldots, b-a-2$\\
* $\partial_{\olz_2} \widetilde{q}_{a,b;0,b-a-1} = q_{a,b-1;0,b-a-1} $\\
* $\partial_{\olz_2} q_{a,b;0,d} = q_{a,b-1;0,d}, \quad d = b-a, \ldots, b-1$\\
* $\partial_{\olz_2} q_{a,b;0,b} = 0$\\[-2mm]

\noindent (ii) if $a>b$\\[-2mm]

\noindent
* $\partial_{z_2} \widetilde{p}_{a,b;c,0} = \widetilde{p}_{a-1,b;c,0}, \quad c = 0, \ldots, a-b-2$\\
* $\partial_{z_2} \widetilde{p}_{a,b;a-b-1,0} = p_{a-1,b;a-b-1,0}$\\
* $\partial_{z_2} p_{a,b;c,0} = p_{a-1,b;c,0}, \quad c = a-b, \ldots, a-1$\\
* $\partial_{z_2} p_{a,b;a,0} = 0$\\
* $\partial_{z_2} q_{a,b;0,d} = q_{a-1,b;0,d}, \quad d = 0, \ldots, b$\\[-2mm]

\noindent
* $\partial_{\olz_2} \widetilde{p}_{a,b;c,0} = \widetilde{p}_{a,b-1;c,0}, \quad c = 0, \ldots, a-b-1$\\
* $\partial_{\olz_2} p_{a,b;a-b,0} = \widetilde{p}_{a,b-1;a-b,0}$\\
* $\partial_{\olz_2} p_{a,b;c,0} = p_{a,b-1;c,0}, \quad c = a-b+1, \ldots, a$\\
* $\partial_{\olz_2} q_{a,b;0,d} = q_{a,b-1;0,d}, \quad d = 0, \ldots, b-1$\\
* $\partial_{\olz_2} q_{a,b;0,b} = 0$\\[-2mm]

\noindent (iii) if $a=b$\\[-2mm]

\noindent
* $\partial_{z_2} p_{a,a;c,0} = p_{a-1,a;c,0}, \quad c = 0, \ldots, a-1$\\
* $\partial_{z_2} p_{a,a;a,0} = 0$\\
* $\partial_{z_2} q_{a,a;0,0} = \widetilde{q}_{a-1,a;0,0}$\\
* $\partial_{z_2} q_{a,a;0,d} = q_{a-1,a;0,d}, \quad d = 1, \ldots, a$\\[-2mm]

\noindent
* $\partial_{\olz_2} p_{a,a;0,0} = \widetilde{p}_{a,a-1;0,0}$\\
* $\partial_{\olz_2} p_{a,a;c,0} = p_{a,a-1;c,0}, \quad c = 1, \ldots, a$\\
* $\partial_{\olz_2} q_{a,a;0,d} = q_{a,a-1;0,d}, \quad d = 0, \ldots, a-1$\\
* $\partial_{\olz_2} q_{a,a;0,a} = 0$\\[-2mm]

\noindent For the derivatives with respect to $z_1$ and $\olz_1$ there holds:\\[-2mm]

\noindent (i) if $a<b$\\[-2mm]

\noindent
*$\partial_{z_1} p_{a,b;0,0} = - \widetilde{q}_{a-1,b;0,0}$\\
*$\partial_{z_1} p_{a,b;c,0} = \phantom{-} p_{a-1,b;c-1,0}, \quad c = 1, \ldots, a$\\
*$\partial_{z_1} \widetilde{q}_{a,b;0,d} = -  \widetilde{q}_{a-1,b;0,d+1}, \quad d = 0, \ldots, b-a-1$\\
*$\partial_{z_1} q_{a,b;0,d} = - q_{a-1,b;0,d+1}, \quad d = b-a, \ldots, b-1$\\
*$\partial_{z_1} q_{a,b;0,b} = 0$\\[-2mm]

\noindent
*$\partial_{\olz_1} p_{a,b;c,0} = - p_{a,b-1;c+1,0}, \quad c = 0, \ldots, a-1$\\
*$\partial_{\olz_1} p_{a,b;a,0} = 0$\\
*$\partial_{\olz_1} \widetilde{q}_{a,b;0,0} = p_{a,b-1;0,0}$\\
*$\partial_{\olz_1} \widetilde{q}_{a,b;0,d} =  \widetilde{q}_{a,b-1;0,d-1}, \quad d = 1, \ldots, b-a-1$\\
*$\partial_{\olz_1} q_{a,b;0,d} =  q_{a,b-1;0,d-1}, \quad d = b-a, \ldots, b$\\[-2mm]

\noindent (ii) if $a>b$\\[-2mm]

\noindent
* $\partial_{z_1} \widetilde{p}_{a,b;0,0} = - q_{a-1,b;0,0}$\\
* $\partial_{z_1} \widetilde{p}_{a,b;c,0} = \widetilde{p}_{a-1,b;c-1,0}, \quad c = 1, \ldots, a-b-1$\\
* $\partial_{z_1} p_{a,b;c,0} = p_{a-1,b;c-1,0} , \quad c = a-b, \ldots, a$\\
* $\partial_{z_1} q_{a,b;0,d} = q_{a-1,b;0,d+1}, \quad d = 0, \ldots, b-1$\\
* $\partial_{z_1} q_{a,b;0,b} = 0$\\[-2mm]

\noindent
* $\partial_{\olz_1} \widetilde{p}_{a,b;c,0} = - \widetilde{p}_{a,b-1;c+1,0}, \quad c = 0, \ldots, a-b-1$\\
* $\partial_{\olz_1} p_{a,b;c,0} = \phantom{-} p_{a,b-1;c+1,0}, \quad c = a-b, \ldots, a-1$\\
* $\partial_{\olz_1} p_{a,b;a,0} = 0$\\
* $\partial_{\olz_1} q_{a,b;0,0} = \widetilde{p}_{a,b-1;0,0}$\\
* $\partial_{\olz_1} q_{a,b;0,d} = q_{a,b-1;0,d-1}, \quad d = 1, \ldots, b$\\[-2mm]

\noindent (iii) if $a=b$\\[-2mm]

\noindent
* $\partial_{z_1} p_{a,a;0,0} = \widetilde{q}_{a-1,a;0,0}$\\
* $\partial_{z_1} p_{a,a;c,0} = \phantom{-} p_{a-1,a;c-1,0}, \quad c = 1, \ldots, a$\\
* $\partial_{z_1} q_{a,a;0,d} = - q_{a-1,a;0,d+1}, \quad d = 0, \ldots, a-1$\\
* $\partial_{z_1} q_{a,a;0,a} = 0$\\

\noindent
* $\partial_{\olz_1} p_{a,a;c,0} = - p_{a,a-1;c+1,0}, \quad c = 0, \ldots, a-1$\\
* $\partial_{\olz_1} p_{a,a;a,0} = 0$\\
* $\partial_{\olz_1} q_{a,a;0,0} = \widetilde{p}_{a,a-1;0,0}$\\
* $\partial_{\olz_1} q_{a,a;0,d} = \phantom{-} q_{a,a-1;0,d-1}, \quad d = 1, \ldots, a$
\end{proposition}

\pf
Follows by direct computation using the properties of the $Q_l^{\alpha, \beta}(\uz, \uzd)$ polynomials established in Proposition 1.
\qed


\section{Conclusion}

The underlying paper may be seen as the conclusive tailpiece in a series of papers on Hermitian Clifford analysis. In \cite{partI} it was shown how Hermitian Clifford analysis arises quite naturally as a special case of standard Clifford analysis by introducing a so--called complex structure, i.e. a special orthogonal matrix, or the corresponding Spin--element in the Clifford algebra, squaring up to $- {\bf 1}$, in this way breaking down the orthogonal invariance of standard Clifford analysis to the unitary one. When considering functions with values in the whole Clifford algebra it is a known fact that the first order system of equations expressing (Hermitian) monogenicity contains redundant information; in \cite{partII}  the conceptual meaning of Hermitian monogenicity was further unraveled by studying possible splittings of the corresponding first order system into independent parts without changing the properties of the solutions, leading to Hermitian monogenic functions with values in spinor space and subspaces thereof; in this way also connections with holomorphic functions of several complex variables were established. A fundamental result which may not be missing in a function theory, is the so--called Fischer decomposition; in \cite{archivum} spaces of homogeneous monogenic polynomials were decomposed into unitary--irreducibles involving homogeneous Hermitean monogenic polynomials. In \cite{howe} the choice of the Hermitian monogenicity equations was fully justified; indeed, constructing the Howe dual for the action of the unitary group on the space of all spinor valued polynomials, the generators of the resulting Lie superalgebra reveal the natural set of equations to be considered in this context, which exactly coincide with the chosen ones.  Next to the Fischer decomposition, a second essential step towards the construction of an orthogonal basis of Hermitean monogenic polynomials, was taken in \cite{CK} by establishing a Cauchy--Kovalevskaya extension theorem for such polynomials. In general the problem of constructing orthogonal bases for spaces of null solutions of a given partial differential operator is quite difficult. However, if the partial differential equation has a sufficiently broad symmetry, i.e., if the group preserving the space of solutions is sufficiently large, the construction is much facilitated by the Gel'fand--Tsetlin approach, which then offers an efficient tool for it. Indeed, the notion of Gel'fand--Tsetlin (GT) basis applies to finite dimensional irreducible modules over a classical Lie algebra; when this module is realized explicitly, say as a subspace of null solutions of an invariant differential operator, then an algorithm for the construction of the GT--basis may be devised. This was accomplished in \cite{rhohds,rhofb} where in a systematic and detailed way the GT--construction of orthogonal bases for spaces of homogeneous Hermitian monogenic polynomials was described with special attention for the Appell property. Finally, in the present paper,
the embedding factors for the step 2 branching of spaces of homogeneous Hermitian monogenic polynomials are explicitly determined, allowing for a new inductive construction of the orthogonal bases for spaces of Hermitian spherical monogenics, with special emphasis on the complex dimension 2 case.


\section*{Acknowledgements}


R. L\'{a}vi\v{c}ka and V. Sou\v{c}ek acknowledge support by the  grant GA CR 201/08/0397.




\begin{thebibliography}{11}
 
\bibitem{bock} S.\ Bock, K.\ Guerlebeck, R.\ L\'{a}vi\v{c}ka, V.\ Sou\v{c}ek, The Gel'fand-Tsetlin bases for spherical monogenics in dimension 3,  arXiv:1010.1615v2 [math.CV], 2010, to appear in \textit{Rev. Mat. Iberoamericana}.

\bibitem{partI}	F.\ Brackx, J.\ Bure\v{s}, H.\ De Schepper, D.\ Eelbode, F.\ Sommen, V.\ Sou\v{c}ek,
Fundaments of Hermitian Clifford analysis -- Part I: Complex structure, {\em Compl. Anal. Oper. Theory}  {\bf 1} (3), 2007, 341--365.

\bibitem{partII} F.\ Brackx, J.\ Bure\v{s}, H.\ De Schepper, D.\ Eelbode, F.\ Sommen, V.\ Sou\v{c}ek, 
Fundaments of Hermitian Clifford analysis -- Part II: Splitting of $h$--monogenic equations, {\em Complex Var. Elliptic Eq.} {\bf 52} (10-11), 2007, 1063--1079.

\bibitem{cama} F.\ Brackx, B.\ De Knock, H.\ De Schepper, F.\ Sommen, 
On Cauchy and Martinelli-Bochner Formulae in Hermitian Clifford Analysis, \textit{Bull. Braz. Math. Soc.} {\bf 40} (3), 2009, 395--416.

\bibitem{mahi} F.\ Brackx, B.\ De Knock, H.\ De Schepper, A matrix Hilbert transform in Hermitian Clifford Analysis, \textit{J. Math. Anal. Appl.} {\bf 344} (2), 2008, 1068--1078.

\bibitem{bds} F.\ Brackx, R.\ Delanghe, F.\ Sommen, \textit{Clifford Analysis}, Pitman Publishers (Boston-London-Melbourne, 1982).

\bibitem{hehe} F.\ Brackx, H.\ De Schepper, N.\ De Schepper, F.\ Sommen,
Hermitian Clifford-Hermite polynomials, {\em Adv. Appl. Clifford Alg.} {\bf 17} (3), 2007, 311--330.

\bibitem{howe} F.\ Brackx, H.\ De Schepper, D.\ Eelbode, V.\ Sou\v{c}ek, 
The Howe Dual Pair in Hermitian Clifford Analysis, {\em Rev. Mat. Iberoamericana} {\bf 26} (2), 2010, 449--479.

\bibitem{CK} F.\ Brackx, H.\ De Schepper, R.\ L\'{a}vi\v{c}ka, V.\ Sou\v{c}ek, 
The Cauchy--Kovalevskaya Extension Theorem in Hermitian Clifford Analysis, {\em J. Math. Anal. Appl.} {\bf 381} (2), 2011, 649--660.

\bibitem{rhovla} F.\ Brackx, H.\ De Schepper, R.\ L\'{a}vi\v{c}ka, V.\ Sou\v{c}ek, 
Gel'fand-Tsetlin procedure for the construction of orthogonal bases in Hermitian Clifford analysis. In: T.E. Simos, G. Psihoyios, Ch. Tsitouras, \textit{Numerical Analysis and Applied Mathematics}, AIP Conference Proceedings 1281, Rhodos, Greece (2010), 1508--1511.

\bibitem{rhoro} F.\ Brackx, H.\ De Schepper, R.\ L\'{a}vi\v{c}ka, V.\ Sou\v{c}ek,
Fischer decompositions of kernels of Hermitian Dirac operators. In: T.E. Simos, G. Psihoyios, Ch. Tsitouras, \textit{Numerical Analysis and Applied Mathematics}, AIP Conference Proceedings 1281, Rhodos, Greece (2010), 1484--1487. 

\bibitem{rhofb} F.\ Brackx, H.\ De Schepper, R.\ L\'{a}vi\v{c}ka, V.\ Sou\v{c}ek, 
Orthogonal bases of Hermitian monogenic polynomials: an explicit construction in complex dimension 2. In: T.E. Simos, G. Psihoyios, Ch. Tsitouras, \textit{Numerical Analysis and Applied Mathematics, AIP Conference Proceedings 1281}, Rhodos, Greece (2010), 1451--1454.

\bibitem{rhohds} F.\ Brackx, H.\ De Schepper, R.\ L\'{a}vi\v{c}ka, V.\ Sou\v{c}ek,
Gel'fand-Tsetlin Bases of Orthogonal Polynomials in Hermitian Clifford Analysis, \textit{Math. Meth. Appl. Sci.} {\bf 34}, 2011, 2167--2180.

\bibitem{toulouse} F.\ Brackx, H.\ De Schepper, F.\ Sommen, The Hermitian Clifford analysis toolbox, \textit{Adv. Appl. Cliff. Alg.} {\bf 18} (3-4), 2008, 451--487.

\bibitem{hermwav} F.\ Brackx, H. De Schepper, F.\ Sommen, A Theoretical Framework for Wavelet Analysis in a Hermitian Clifford Setting, {\em Communications on Pure and Applied Analysis} {\bf 6} (3), 2007, 549--567.

\bibitem{archivum} F.\ Brackx, H.\ De Schepper, V.\ Sou\v{c}ek,
Fischer  Decompositions in Euclidean and Hermitian Clifford Analysis,
{\em Archivum Mathematicum (Brno)} {\bf 46} (5), 2010, 301--321.

\bibitem{cacao} I.\ Ca\c{c}\~{a}o,
\textit{Constructive approximation by monogenic polynomials}, PhD thesis, University of Aveiro (2004).

\bibitem{cagubo} I.\ Ca\c{c}\~{a}o, K.\ G\"{u}rlebeck, S. Bock,
Complete orthonormal systems of spherical monogenics -- a constructive approach. In: L.\ H.\ Son. W. Tutschke, S. Jain (eds.), \textit{Methods of Complex and Clifford Analysis}, Proceedings of ICAM, SAS Int. Pub. (Hanoi, 2004).

\bibitem{struppa} F.\ Colombo, I.\ Sabadini, F.\ Sommen, D.\ C.\ Struppa, \textit{Analysis of Dirac Systems and Computational Algebra}, Birkh\"{a}user (Boston, 2004).

\bibitem{damdeef} A.\ Damiano and D.\ Eelbode, Invariant Operators Between Spaces of h-Monogenic Polynomials, \textit{Adv. Appl. Cliff. Alg.} {\bf 19} (2), 2009, 237--251.

\bibitem{dss} R.\ Delanghe, F.\ Sommen, V.\ Sou\v{c}ek, 
\textit{Clifford algebra and spinor-valued functions -- A function theory for the Dirac operator}, Kluwer Academic Publishers (Dordrecht, 1992).

\bibitem{eel} D.\ Eelbode, 
Stirling numbers and Spin--Euler polynomials, {\em Exp. Math.} {\bf 16} (1), 2007, 55--66.

\bibitem{eel2} D.\ Eelbode, 
Irreducible $\gsl(m)$--modules of Hermitian monogenics, \textit{Complex Var. Elliptic Equ.} {\bf 53} (10), 2008, 975--987.

\bibitem{eel3} D.\ Eelbode, Fu Li He,
Taylor series in Hermitian Clifford Analysis, \textit{Complex Analysis and Operator Theory}, DOI: 10.1007/s11785-009-0036-y.

\bibitem{gilbert} J.\ Gilbert, M.\ Murray, 
\textit{Clifford Algebra and Dirac Operators in Harmonic Analysis}, Cambridge University Press (Cambridge, 1991).

\bibitem{ghs} K.\ G\"{u}rlebeck, K.\ Habetha, W.\ Spr\"{o}\ss ig, 
\textit{Holomorphic functions in the plane and $n$-dimensional space. Translated from the 2006 German original. With 1 CD-ROM (Windows and UNIX)}, Birkh\"{a}user Verlag (Basel, 2008).

\bibitem{guerleb} K.\ G\"urlebeck, W.\ Spr\"o\ss ig, 
\textit{Quaternionic and Clifford Calculus for Physicists and Engineers}, J.\ Wiley \& Sons (Chichester, 1997).

\bibitem{lav} R.\ L\'{a}vi\v{c}ka,
Canonical bases for $\mathfrak{s}\mathfrak{l}(2,\mC)$-modules of spherical monogenics in dimension 3, \textit{Archivum Mathematicum} {\bf 46} (5) 2010, 339--349. 

\bibitem{lav1} R.\ L\'{a}vi\v{c}ka, Complete orthogonal Appell systems for spherical monogenics, \textit{Compl. Anal. Oper. Theory} {\bf 6}  (2) 2012, 477--489.

\bibitem{lav2} R.\ L\'{a}vi\v{c}ka, Orthogonal Appell bases for Hodge-de Rham systems in Euclidean spaces, arXiv:1111.0974v1 [math.CV], 2011, to appear in \textit{Adv. Appl. Clifford Alg.} 

\bibitem{lasova} R.\ L\'{a}vi\v{c}ka, V.\ Sou\v{c}ek, P.\ Van Lancker,
Orthogonal basis for spherical monogenics by step two branching, \textit{Ann. Glob. Anal. Geom.} {\bf 41} (2), 2012, 161-186.

\bibitem{molev} A.\ I.\ Molev,
Gel'fand--Tsetlin bases for classical Lie algebras. In: M.\ Hazewinkel (ed.), \textit{Handbook of Algebra, vol.4}, Elsevier, (2006), 109--170.

\bibitem{port} I.\ Porteous,
\textit{Clifford Algebras and the Classical groups}, Cambridge University Press (Cambridge, 1995).

\bibitem{rocha} R.\ Rocha-Chavez, M.\ Shapiro, F.\ Sommen, 
\textit{Integral theorems for functions and differential forms in $\mC_m$}, Research Notes in Math. 428, Chapman\&Hall / CRC (New York, 2002).

\bibitem{sabadini} I.\ Sabadini, F.\ Sommen, 
Hermitian Clifford analysis and resolutions, {\em Math. Meth. Appl. Sci.} {\bf 25}(16-18), 2002, 1395--1414.

\bibitem{mabo} F.\ Sommen, D.\ Pe\~{n}a Pe\~{n}a, A Martinelli-Bochner formula for the Hermitian Dirac equation, {\em Math. Meth. Appl. Sci.}  {\bf 30} (9), 2007, 1049--1055.

\bibitem{pvl} P.\ Van Lancker,
Spherical Monogenics: An Algebraic Approach, {\em Adv. Appl. Clifford Alg.} {\bf 19}, 2009, 467--496.

\end{thebibliography}
\end{document}